\documentclass[a4paper,12pt]{article}
\usepackage{amssymb,amsmath,amsthm,}
\newtheorem{proposition}{Proposition}[section]
\newtheorem{theorem}{Theorem}[section]
\newtheorem{lemma}{Lemma}[section]
\newtheorem{corollary}{Corollary}
\newtheorem{remark}{Remark}[section]
\begin{document}
\title {Asymptotics of Eigenvalues and Eigenfunctions for the
Laplace Operator in a Domain with Oscillating Boundary.\\
Multiple Eigenvalue Case.$\!\!$ \footnote{The work of the second
author was partially supported by RFBR (06-01-00441) and by the
program ``Leading Scientific Schools'' (HIII-2538.2006.1). The
work of the third author was partially supported by RFBR
(06-01-00138).}}
\author{Youcef Amirat$^\flat$, Gregory A. Chechkin$^\natural$, Rustem R. Gadyl'shin$^\sharp$.}
\date{}
\maketitle
\bigskip
\begin{center}\vskip-15pt
{\footnotesize $^\flat$ Laboratoire de Math\'ematiques Appliqu\'ees\\
CNRS UMR 6620\\
Universit\'e Blaise Pascal\\
63177 Aubi\`ere cedex, France\\
{\tt amirat@math.univ-bpclermont.fr}
\bigskip
$^\natural$ Department of Differential Equations\\
Faculty of Mechanics and Mathematics\\
Moscow State University\\
Moscow 119992, Russia\\
{\tt chechkin{@}mech.math.msu.su}
\bigskip
$^\sharp$ 
Department of Mathematical Analysis\\ Faculty of Physics and
Mathematics\\ Bashkir State Pedagogical University\\ Ufa 450000,
Russia\\ {\tt gadylshin{@}yandex.ru,\, gadylshin{@}bspu.ru}}
\end{center}

\newpage\vbox{}
\centerline{Abstract}

\bigskip
\noindent We study the asymptotic behavior of the solutions of a
spectral problem for the Lap\-la\-ci\-an in a domain with
rapidly oscillating boundary. We consider the case where the
eigenvalue of the limit problem is multiple. We construct the
leading terms of the asymptotic expansions for the eigenelements
and verify the asymptotics.
\vskip 50pt { \centerline{R\'esum\'e}
\bigskip
\noindent Nous \'etudions le comportement asymptotique des
solutions d'un probl\`eme spectral associ\'e \`a l'op\'erateur
de Laplace dans un domaine \`a fronti\`ere oscillante. Nous
consid\'erons le cas o\`u la valeur propre du probl\`eme limite
est multiple. Nous construisons les termes principaux des
d\'eveloppements asymptotiques des \'el\'ements propres et nous
donnons une justification rigoureuse des ap\-pro\-xi\-ma\-tions
asymptotiques.}
\vskip 10.cm
\def \trait (#1) (#2) (#3){\vrule width #1pt height #2pt depth #3pt}

\noindent{{\trait (156) (0.1) (0.1)}\hfill}

{\footnotesize {\it Keywords.} spectral problems, oscillating
boundary, asymptotic expansions


}

\newpage
\section* {Introduction.}
Boundary-value problems involving oscillating boundaries or
interfaces appear in many fields of natural sciences and
engineering, such as the scattering of acoustic and
electro--magnetic waves on small periodic obstacles (for
instance, whispering gallery effects \cite{Pe} and scanning of
the surface of oceans from the outer space \cite{BMS}), the
vibrations of strongly inhomogeneous elastic bodies \cite{BPC}
\cite{CC1}, the friction of details in complex engineering
structures \cite{B}, the flows over rough walls \cite{APV}, or
behavior of coupled fluid-solid periodic structures (structures
having soft and hard phases \cite{BLS}). Recent years many other
ma\-the\-ma\-ti\-cal works (purely theoretical as well as
applied) were devoted to asymptotic analysis of these problems,
see for instance, \cite{AS}, \cite{BV},
\cite{CCHE1}--\cite{FHL}, \cite{Ga1}, \cite{JwMa}, \cite{KPV},
\cite{MK}--\cite{NK}, \cite{SP}.

In the paper \cite{BV} the authors considered spectral problems
for a general $2m$-order elliptic operator in a domain of a
special type with partially oscillating boundary with the
Dirichlet type of boundary conditions on the oscillating part of
the boundary. They proved the convergence theorem for the
eigenvalues and eigenfunctions. Also it should be noted that
similar convergence results were given in \cite{O}, as
application of the method for the approximation of eigenvalues
and eigenvectors of self-adjoint operators.

In \cite{AmCheGa} we considered a spectral problem for the
Laplace operator in a bounded domain with the boundary which
part, depending on a small pa\-ra\-me\-ter $\varepsilon$, is
rapidly oscillating. The authors assumed that the frequency and
the amplitude of oscillations of the boundary are of the same
order $\varepsilon$. The case of simple eigenvalue of the limit
problem is studied:  the authors constructed the leading terms
of the asymptotic expansions for the eigenelements and verified
the asymptotics.

In this paper we deal with the same spectral problem in the case
when the eigenvalue of the limit problem is multiple. Our aim is
to construct accurate asymptotic approximations, as $\varepsilon
\to 0$, of the eigenvalues and corresponding eigenfunctions. We
use the method of matching of asymptotic expansions
(see~\cite{I0},~\cite{I1} and~\cite{I3}) to construct the leading
terms of the asymptotic expansions for the eigenelements. Then we
prove the asymptotic estimates of the difference between the
solutions of the original problem and the ap\-pro\-xi\-mate
asymptotic expansions (see also papers \cite{CUMN} and \cite{C}).

The case of the domain with totally oscillating boundary is
considered in \cite{N}. For such a domain the author constructed
two terms asymptotics of the eigenvalues. Neumann boundary-value
problems were considered in \cite{LS} and also in \cite{NM2},
\cite{NO}.

The outline of the paper is as follows: in Section 1 we
introduce the notations, set the problem, give preliminary
propositions and statements of the main results. In Section 2 we
derive the formal asymptotics for the eigenelements, in Section
3 we give a rigorous justification of the asymptotics, and in
Appendix we prove two auxiliary Propositions.

\section {Setting of the problem,
preliminary pro\-po\-si\-tions and statements of the main
results}

Let $\Omega$ be a bounded domain in $\mathbb{R}^2$, located in the
upper half space. We assume the boundary $\partial\Omega$ to be
piecewise smooth, consisting of the parts:
$\partial\Omega=\Gamma_0\cup\Gamma_1\cup\Gamma_2\cup\Gamma_3$,
where $\Gamma_0$ is the segment $(-\frac{1}{2},\frac{1}{2})$ on
the abscissa axis, $\Gamma_2$ and $\Gamma_3$ belong to the
straight lines  $x_1=-\frac{1}{2}$ and $x_1=\frac{1}{2}$,
respectively. Let $\varepsilon=\frac{1}{2\mathcal{N}+1}$ be a
small parameter, where $\mathcal{N}$ is a large positive number.
Given a smooth negative $1$-periodic in\, $\xi_1$\, even function
\, $F(\xi_1)$, such that $F'(\xi_1)= 0$ for $\xi_1=\pm\frac{1}{2}$
and $\xi_1=0$, we set $$ \Pi_\varepsilon = \{ x \in \mathbb{R}^2 \
: \ (x_1,0) \in \Gamma_0, \ \varepsilon
F\left(\frac{x_1}{\varepsilon}\right) < x_2 \leq 0 \} $$ and then
we denote  $$\Omega^\varepsilon = \Omega \cup \Pi_\varepsilon.$$
Thus, the boundary of $\Omega^\varepsilon$ consists of four parts:
$\partial\Omega^\varepsilon=\Gamma_\varepsilon\cup\Gamma_1\cup
\Gamma_{2,\varepsilon}\cup\Gamma_{3,\varepsilon}$, where $$
\Gamma_\varepsilon=\{ x \in \mathbb{R}^2 \ : \ (x_1,0) \in
\Gamma_0, \ x_2=\varepsilon
F\left(\frac{x_1}{\varepsilon}\right)\}, $$ $$
\Gamma_{2,\varepsilon}=\Gamma_2\cup\{x \in \mathbb{R}^2 \ : \
x_1=-\frac{1}{2}, \ \varepsilon
F\left(-\frac{1}{2\varepsilon}\right)\leq x_2\leq0\}, $$ $$
\Gamma_{3,\varepsilon}=\Gamma_3\cup\{x \in \mathbb{R}^2 \ : \
x_1=\frac{1}{2}, \ \varepsilon
F\left(\frac{1}{2\varepsilon}\right)\leq x_2\leq0\}. $$ Denote by
$\Gamma=\{\xi\in \mathbb{R}^2\ :\ -\frac{1}{2}<\xi_1<\frac{1}{2},
\ \xi_2=F(\xi_1)\}$ in $\xi=\frac{x}{\varepsilon}$ variables, and
let $\Pi=\{\xi\in \mathbb{R}^2\ :\ -\frac{1}{2}<\xi_1<\frac{1}{2},
\ \xi_2>F(\xi_1)\}$ be a semi-infinite strip.

Denote 
by $\nu$ the outward unit normal vector. The following statement
is proved in \cite{AmCheGa}.

\begin{lemma}\label{C}
Assume that the multiplicity of the eigenvalue $\lambda_0$ of
Problem
\begin{equation}\label{2}
\left\{\begin{array}{l} -\Delta u_0=\lambda_0 u_0\quad\hbox{in}
\ \Omega,\\ u_0=0\quad\hbox{on} \ \Gamma_0,\\ \frac{\partial
u_0}{\partial \nu}=0\quad\hbox{on} \
\Gamma_1\cup\Gamma_2\cup\Gamma_3.
\end{array}\right.
\end{equation}
is equal to $p$. Then there are $p$ eigenvalues of Problem
\begin{equation}\label{1}
\left\{\begin{array}{l} -\Delta
u_\varepsilon=\lambda_\varepsilon u_\varepsilon\quad\hbox{in} \
\Omega^\varepsilon,\\ u_\varepsilon=0\quad\hbox{on} \
\Gamma_\varepsilon,\\ \frac{\partial u_\varepsilon}{\partial
\nu}=0\quad\hbox{on} \
\Gamma_1\cup\Gamma_{2,\varepsilon}\cup\Gamma_{3,\varepsilon}.
\end{array}\right.
\end{equation}
(with multiplicities taken into account) converging to
$\lambda_0$, as $\varepsilon\to 0$.
\end{lemma}

In \cite{AmCheGa} we considered the case where $\lambda_0$ is
simple:  we constructed the leading terms of the asymptotic
expansions for the eigenelements and verified the asymptotics.
Here we deal with the case where  $\lambda_0$ is multiple.

Here and throughout we assume, without loss of generality,  that
the multiplicity of $\lambda_0$ equals two. Let then $u_0^{(l)}$
$(l=1, 2)$ be the basis of the eigensubspace corresponding to
$\lambda_0$, formed by eigenfunctions of Problem (\ref{2}),
orthonormalized in $L_2(\Omega)$:
\begin{equation*}\label{2p1}
\left\{\begin{array}{l} -\Delta u_0^{(l)}=\lambda_0
u_0^{(l)}\quad\hbox{in} \ \Omega,\\ u_0^{(l)}=0\quad\hbox{on} \
\Gamma_0,\\ \frac{\partial u_0^{(l)}}{\partial
\nu}=0\quad\hbox{on} \ \Gamma_1\cup\Gamma_2\cup\Gamma_3,
\end{array}\right.
\end{equation*}

\begin{equation}\label{2p2}
\int \limits_{\Omega} (u_0^{(l)})^2 \, dx=1, \quad
\int_{\Omega} u_0^{(1)}\, u_0^{(2)}  \, dx=0, \quad l=1, 2.
\end{equation}

It is easy to see that the eigenvalues can be chosen to satisfy
an additional orthogonality condition on $\Gamma_0$
\begin{equation}\label{6.26}
\int\limits_{\Gamma_0}\frac{\partial u_0^{(1)}}{\partial \nu}
\frac{\partial u_0^{(2)}}{\partial \nu}\ ds=0.
\end{equation}
Note that the similar orthogonality condition on the boundary of
the type~(\ref{6.26}) was used in~\cite{Bor} and~\cite{C}. In
addition for simplicity we assume that
\begin{equation}\label{nerav}
\int\limits_{-\frac{1}{2}}^{\frac{1}{2}}\left(\frac{\partial
u_0^{(1)}}{\partial x_2}\right)^2 d x_1\neq
\int\limits_{-\frac{1}{2}}^{\frac{1}{2}}\left(\frac{\partial
u_0^{(2)}}{\partial x_2}\right)^2 d x_1.
\end{equation}
Due to Lemma \ref{C}, there are two eigenvalues of Problem
(\ref{1}), denoted $\lambda_\varepsilon^{(1)}$ and
$\lambda_\varepsilon^{(2)}$, converging to $\lambda_0$, as
$\varepsilon\to 0$. Throughout we denote by
$u_\varepsilon^{(l)}$ $(l=1, 2)$ the corresponding
eigenfunctions, orthonormalized in $L_2(\Omega_\varepsilon)$. We
then have

\begin{equation}\label{2p3}
\left\{\begin{array}{l} -\Delta
u_\varepsilon^{(l)}=\lambda_\varepsilon^{(l)}
u_\varepsilon^{(l)}\quad\hbox{in} \ \Omega^\varepsilon,\\
u_\varepsilon^{(l)}=0\quad\hbox{on} \ \Gamma_\varepsilon,\\
\frac{\partial u_\varepsilon^{(l)}}{\partial
\nu}=0\quad\hbox{on} \
\Gamma_1\cup\Gamma_{2,\varepsilon}\cup\Gamma_{3,\varepsilon}.
\end{array}\right.
\end{equation}

\begin{equation}\label{2p4}
\int \limits_{\Omega_\varepsilon} (u_\varepsilon^{(l)})^2 \,
dx=1, \quad \int \limits_{\Omega_\varepsilon}
u_\varepsilon^{(1)}\, u_\varepsilon^{(2)}\, dx=0, \quad l=1, 2.
\end{equation}

Our aim is the construction of accurate asymptotic
approximations, as $\varepsilon \to 0$, for the eigenvalues
$\lambda_\varepsilon^{(1)}$ and $\lambda_\varepsilon^{(2)}$ and
for the corresponding eigenfunctions.

It is proved in \cite{AmCheGa} that Problem
\begin{equation}\label{6'}
\left\{\begin{array}{l} \Delta_\xi X=0\quad\hbox{in} \ \Pi,\\
X=0\quad\hbox{on} \ \Gamma,\qquad \frac{\partial X}{\partial
\xi_1}=0\quad\hbox{as} \ \xi_1=\pm\frac{1}{2},\\ X\sim \xi_2
\quad\hbox{as}\ \xi_2\to+\infty.
\end{array}\right.
\end{equation}
has a solution  with the asymptotics
\begin{equation}\label{7'}
X(\xi)=\xi_2+C(F) \quad\hbox{as}\ \xi_2\to+\infty,
\end{equation}
up to exponentially small terms, where $C(F)$ is a positive
constant depending on the function $F$. Note that, due to the
evenness of the function $F$, the solution $X$ is even in
$\xi_1$ and can be extended to a $1$-periodic function in
$\xi_1$. Later on we use the same notation $X$ for the
extension.

Our main goal is to prove the following statement.

\begin{theorem}\label{t2}
Assume that the multiplicity of $\lambda_0$ of Problem
{\rm(\ref{2})} equals two, the associated eigenfunctions
$u_0^{(l)}$ ($l=1, 2$) satisfy the conditions
{\rm(\ref{2p2})--(\ref{nerav})}. Then eigenvalues
$\lambda_\varepsilon^{(l)}$  of Problem {\rm(\ref{1})}, converging
to $\lambda_0$ as $\varepsilon \to 0$, and the associated
eigenfunctions $u_\varepsilon^{(l)}$ orthonormalized in
$L_2(\Omega_\varepsilon)$ have the following asymptotics:
\begin{align}\label{est1}
\lambda_\varepsilon^{(l)}=&\lambda_0+\varepsilon
\lambda_1^{(l)}+o\left(\varepsilon^{\frac{5}{4}-\sigma}\right)\qquad\hbox{for
any $\sigma>0$},\\ \lambda_1^{(l)}=& -C(F)
\int\limits_{\Gamma_0}\left(\frac{\partial u_0^{(l)}}{\partial
\nu}\right)^2\ ds, \label{Est1}\\
&\|u^{(l)}_\varepsilon-u^{(l)}_0\|_{H^1(\Omega)} +
\|u^{(l)}_\varepsilon\|_{H^1(\Omega^\varepsilon\backslash\overline{\Omega})}
=o(1).\label{Est2}
\end{align}
\end{theorem}

\begin{remark}\label{rem1}
In the next section we construct four terms asymptotics of the
eigenelements of Problem (\ref{2p3}). Moreover given algorithm
allows to construct (see Remark~\ref{rem2}) and to justify (see
Remark~\ref{rem3}) the complete asymptotic expansions of the
eigenvalues and eigenfunctions.
\end{remark}

\section {Formal construction of the asymptotics.}\label{s1}
In this section we formally construct the asymptotics of the
eigenvalues $\lambda_\varepsilon^{(l)}$ $(l=1, 2)$ converging to
$\lambda_0$ as $\varepsilon \to 0$, and the asymptotics for
corresponding eigenfunctions $u_\varepsilon^{(l)}$. We use the
method of matching of asymptotic expansions (see
\cite{I0}--\cite{I3}, \cite{G1}--\cite{G3} and
also~\cite{AmCheGa}).
We construct the asymptotics outside a small neighborhood of
$\Gamma_0$ ({\it external} expansion) in the form:
\begin{equation}\label{uu}
u_\varepsilon^{(l)}(x)=u_0^{(l)}(x)+\varepsilon
u_1^{(l)}(x)+\varepsilon^2u_2^{(l)}(x)+\varepsilon^3u_3^{(l)}(x)+
\sum\limits_{i=4}^\infty\varepsilon^iu_i^{(l)}(x),
\end{equation}
the series for the eigenvalues as follows:
\begin{equation}\label{-1}
\lambda_\varepsilon^{(l)}=\lambda_0+\varepsilon\lambda_1^{(l)}+
\varepsilon^2\lambda_2^{(l)}+\varepsilon^3\lambda_3^{(l)}+
\sum\limits_{i=4}^\infty\varepsilon^i\lambda_i^{(l)}
\end{equation}
and the expansion in a small neighborhood of $\Gamma_0$ ({\it
inner} expansion) in the form:
\begin{equation}\label{ve}
u_\varepsilon^{(l)}(x)=\varepsilon v_1^{(l)}(\xi; x_1)+
\varepsilon^2 v_2^{(l)}(\xi; x_1)+\varepsilon^3v_3^{(l)}(\xi;
x_1)+\sum\limits_{i=4}^\infty\varepsilon^iv_i^{(l)}(\xi;x_1),\,
\end{equation}
where $\xi=\frac{x}{\varepsilon}$. Substituting~(\ref{uu})
and~(\ref{-1}) in~(\ref{2p3}) we deduce that the coefficients
of~(\ref{uu}) are to satisfy the following equations and
boundary conditions
\begin{equation}\label{90}
\left\{\begin{array}{l} -\Delta u_1^{(l)}=\lambda_0
u_1^{(l)}+\lambda_1^{(l)} u_0^{(l)}\quad\hbox{in} \ \Omega,\\
\frac{\partial u_1^{(l)}}{\partial \nu}=0\quad\hbox{on} \
\Gamma_1\cup\Gamma_2\cup\Gamma_3,
\end{array}\right.
\end{equation}
\begin{equation}\label{91}
\left\{\begin{array}{l} -\Delta u_2^{(l)}=\lambda_0
u_2^{(l)}+\lambda_1^{(l)} u_1^{(l)}+\lambda_2^{(l)}u_0^{(l)}\quad\hbox{in} \ \Omega,\\
\frac{\partial u_2^{(l)}}{\partial \nu}=0\quad\hbox{on} \
\Gamma_1\cup\Gamma_2\cup\Gamma_3,
\end{array}\right.
\end{equation}
\begin{equation}\label{92}
\left\{\begin{array}{l} -\Delta u_3^{(l)}=\lambda_0
u_3^{(l)}+\lambda_1^{(l)} u_2^{(l)}+\lambda_2^{(l)}u_1^{(l)}+
\lambda_3^{(l)}u_0^{(l)}\quad\hbox{in} \ \Omega,\\
\frac{\partial u_3^{(l)}}{\partial \nu}=0\quad\hbox{on} \
\Gamma_1\cup\Gamma_2\cup\Gamma_3.
\end{array}\right.
\end{equation}
To complete the problems we add boundary conditions on $\Gamma_0:$
\begin{equation}\label{93}
u_i^{(l)}= \alpha_{i0}^{(l)}\quad\hbox{on} \ \Gamma_0,\qquad
i=1,2,\dots,
\end{equation}
where $\alpha_{i0}^{(l)}(x_1)$ are unknown functions, satisfying
the conditions:
\begin{equation}\label{699}
\frac{d^{2k+1}\alpha_{i0}^{(l)}}{dx_1^{2k+1}}\bigg|_{x_1
=\pm\frac{1}{2}}=0,\qquad k=0,1,\dots
\end{equation}
We shall find these functions later. The condition (\ref{699}) is
necessary for solvability of recurrent system of boundary value
problems (\ref{90})--(\ref{93}) in $C^\infty(\overline{\Omega})$.
Moreover such solutions do exist if these problems are solvable in
$H^1(\Omega)$ and in addition due to boundary value problems the
following formulae
\begin{equation}\label{700}
\frac{d^{2k+1}\alpha_{ij}^{(l)}}{dx_1^{2k+1}}\bigg|_{x_1
=\pm\frac{1}{2}}=0,\qquad k=0,1,\dots
\end{equation}
are true, where
\begin{equation}\label{lll}
\alpha_{ij}^{(l)}(x_1)=\frac{1}{j!}\frac{\partial^j
u_i^{(l)}}{\partial x_2^j}\bigg|_{x_2=0}.
\end{equation}
Also it should be noted that due to Problem (\ref{2}) the
following formula
\begin{equation}\label{701}
\alpha_{02}^{(l)}(x_1)\equiv 0
\end{equation}
holds.

Note that, if $\mathcal{F}\in H^1(\Omega)$ and $\alpha\in
H^{1/2}(\Gamma_0)$, then  for solvability in $H^1(\Omega)$ of
the boundary value problem
\begin{equation}\label{90!}
\left\{\begin{array}{l} -\Delta u=\lambda_0
u+\mathcal{F}\quad\hbox{in} \ \Omega,\\
u= \alpha\quad\hbox{on} \ \Gamma_0,\\ \frac{\partial u}{\partial
\nu}=0\quad\hbox{on} \ \Gamma_1\cup\Gamma_2\cup\Gamma_3.
\end{array}\right.
\end{equation}
it is necessary and sufficient to have two identities:
\begin{equation}\label{solva}
\int\limits_{\Omega}\mathcal{F} u_0^{(l)}\
dx=\int\limits_{\Gamma_0}\alpha \frac{\partial
u_0^{(l)}}{\partial \nu}\  ds,\qquad l=1,2.
\end{equation}

In analogues way we obtain boundary value problems for
$v_i^{(l)}$.
\begin{remark}\label{r1}
Further we construct the coefficients of the internal expansion
{\rm (\ref{ve})} in the form of $1$-periodic functions in $\xi_1$.
\end{remark}
In $(\xi,x_1)$ variables the Laplacian and the normal derivative
operator have the form
\begin{equation}\label{5}
\Delta=\varepsilon^{-2} \Delta_\xi
+2\varepsilon^{-1}\frac{\partial^2}{\partial
x_1\partial\xi_1}+\frac{\partial^2}{\partial x_1^2},
\end{equation}
\begin{equation}\label{5+1}
\frac{\partial}{\partial \nu}=\varepsilon^{-1}\frac{\partial
}{\partial \xi_1}+\frac{\partial}{\partial x_1}\quad\hbox{on} \
\Gamma_3,\quad \frac{\partial}{\partial
\nu}=-\varepsilon^{-1}\frac{\partial}{\partial
\xi_1}-\frac{\partial}{\partial x_1}\quad\hbox{on} \ \Gamma_2.
\end{equation}
Substituting~(\ref{-1}),~(\ref{ve}) and keeping in
mind~(\ref{5}),~(\ref{5+1}) and Remark~\ref{r1}, we get the
following equations and boundary conditions for $v_i^{(l)}:$
\begin{equation}\label{60}
\left\{\begin{array}{l} \Delta_\xi v_1^{(l)}=0\quad\hbox{in} \ \Pi,\\
v_1^{(l)}=0\quad\hbox{on} \ \Gamma,\\
\frac{\partial v_1^{(l)}}{\partial \xi_1}=0\quad\hbox{as} \
\xi_1=\pm\frac{1}{2},\ x_1=\pm \frac{1}{2}
.
\end{array}\right.
\end{equation}
\begin{equation}\label{61'}
\left\{\begin{array}{l} -\Delta_\xi v_2^{(l)}= 2\frac{\partial^2
v_1^{(l)}}{\partial x_1\partial\xi_1}\quad\hbox{in} \ \Pi,\\
v_2^{(l)}=0\quad\hbox{on} \ \Gamma,\\
\frac{\partial v_2^{(l)}}{\partial \xi_1}=-\frac{\partial
v_1^{(l)}}{\partial x_1}\quad\hbox{as} \ \xi_1=\pm\frac{1}{2},\
x_1=\pm \frac{1}{2}
.
\end{array}\right.
\end{equation}
\begin{equation}\label{62'}
\left\{\begin{array}{l} -\Delta_\xi v_3^{(l)}= 2\frac{\partial^2
v_2^{(l)}}{\partial x_1\partial\xi_1}+
\frac{\partial^2 v_1^{(l)}}{\partial x_1^2}+\lambda_0 v_1^{(l)}\quad\hbox{in} \ \Pi,\\
v_3^{(l)}=0\quad\hbox{on} \ \Gamma,\\
\frac{\partial v_3^{(l)}}{\partial \xi_1}=-\frac{\partial
v_2^{(l)}}{\partial x_1}\quad\hbox{as} \ \xi_1=\pm\frac{1}{2},\
x_1=\pm \frac{1}{2}
.
\end{array}\right.
\end{equation}
To complete the problems we need to add the conditions at
infinity (as $\xi_2\to+\infty$). These conditions we shall get
from matching of asymptotic expansions. Rewriting the
asymptotics of leading terms of~(\ref{uu}) as $x_2\to0$ in the
variables $\xi=\frac{x}{\varepsilon}$, bearing in mind
(\ref{701}), we deduce
\begin{equation}\label{900}
\sum\limits_{i=0}^{3}\varepsilon^i
u_i^{(l)}(x)=\sum\limits_{i=1}^{3}\varepsilon^i V_i^{(l)}(\xi;
x_1)+O\left(\varepsilon^4 (\xi_2^4+\xi_2)\right)\qquad\hbox{as}\ \
x_2=\varepsilon\xi_2\to0,
\end{equation}
where
\begin{equation}\label{600}
V_1^{(l)}= \alpha_{01}^{(l)}(x_1)\xi_2+\alpha^{(l)}_{10}(x_1),
\end{equation}
\begin{equation}\label{601}
V_2^{(l)}= \alpha_{11}^{(l)}(x_1) \xi_2+\alpha^{(l)}_{20}(x_1),
\end{equation}
\begin{equation}\label{602}
V_3^{(l)}= \alpha_{03}^{(l)}(x_1)\xi_2^3+ \alpha_{12}^{(l)}(x_1)
\xi_2^2+\alpha_{21}^{(l)}(x_1) \xi_2+\alpha^{(l)}_{30}(x_1).
\end{equation}

We must find such $\lambda_i^{(l)}$ and $\alpha^{(l)}_{i0}(x_1)$
that:

$\bullet$\qquad Problems~(\ref{90}) --~(\ref{93}) are to be
soluble,

$\bullet$\qquad Problems~(\ref{60}) --~(\ref{62'}) are to be
soluble with solutions having the asymptotics
\begin{equation}\label{800}
v_i^{(l)}\sim V_i^{(l)} \quad\hbox{as}\ \xi_2\to+\infty
\end{equation}
up to exponentially small terms.

Let us define $\alpha^{(l)}_{10}(x_1)$ and $v_1^{(l)}(\xi; x_1)$.
Obviously the function
\begin{equation}\label{1!}
v_1^{(l)}(\xi;x_1)=\alpha_{01}^{(l)}(x_1) X(\xi)
\end{equation}
due to (\ref{6'}) is the $1$-periodic solution of Problem
(\ref{60})  and due to (\ref{7'}) has the asymptotics
\begin{equation}\label{7}
v_1^{(l)}(\xi;x_1)= \alpha_{01}^{(l)}(x_1) (\xi_2+C(F))
\quad\hbox{as}\ \xi_2\to+\infty,
\end{equation}
up to exponentially small terms. Thus, letting
\begin{equation}\label{702}
\alpha^{(l)}_{10}(x_1)=C(F) \alpha_{01}^{(l)}(x_1),
\end{equation}
we obtain that $v_1^{(l)}$ defined by~(\ref{1!}) satisfies
(\ref{800}), (\ref{600}). Finally, we constructed
$\alpha^{(l)}_{10}$ and $v_1^{(l)}.$

Note that due to (\ref{700}), (\ref{1!})
\begin{equation}\label{5+2}
\frac{\partial v_1^{(l)}}{\partial x_1}
\left(\xi;x_1\right)=0\quad\hbox{as} \ x_1=\pm\frac{1}{2}.
\end{equation}
Hence Problem (\ref{61'}) has the form
\begin{equation}\label{61}
\left\{\begin{array}{l} -\Delta_\xi v_2^{(l)}= 2\frac{\partial^2
v_1^{(l)}}{\partial x_1\partial\xi_1}\quad\hbox{in} \ \Pi,\\
v_2^{(l)}=0\quad\hbox{on} \ \Gamma,\\
\frac{\partial v_2^{(l)}}{\partial \xi_1}=0\quad\hbox{as} \
\xi_1=\pm\frac{1}{2},\ x_1=\pm \frac{1}{2}
.
\end{array}\right.
\end{equation}
and by (\ref{5+1}), (\ref{5+2}) and the boundary condition from
(\ref{6'}), we have
\begin{equation}\label{5+3}
\frac{\partial v_1^{(l)}}{\partial \nu}
\left(\frac{x}{\varepsilon};x_1\right)=0\quad\hbox{on} \
\Gamma_{2,\varepsilon}\cup\Gamma_{3,\varepsilon}.
\end{equation}

Let us define $\lambda_1^{(l)}$ and $u_1^{(l)}(x).$ The constant
$\lambda_1^{(l)}$ can be defined from the solvability condition
of Problem (\ref{90}),~(\ref{93}), which has the same form as
Problem (\ref{90!}). From (\ref{solva}), (\ref{702}),
(\ref{2p2}) and (\ref{6.26}) we deduce that the sufficient
solvability condition of Problem (\ref{90}),~(\ref{93}) is
\begin{equation*}\label{11}
\lambda_1^{(l)}=-C(F) \int\limits_{-\frac{1}{2}}^{\frac{1}{2}}
(\alpha_{01}^{(l)})^2(x_1)\ dx_1
\end{equation*}
or (\ref{Est1}) (taking in account (\ref{lll})). We choose
$u^{(l)}_1$ in the form:
\begin{equation}\label{ttr}
u^{(l)}_1=\widetilde u^{(l)}_1+\kappa_1^{(l)} u^{(l^*)}_0,
\end{equation}
where $$ \int\limits_{\Omega}\widetilde u^{(l)}_1(x)
u^{(k)}_0(x)\ dx=0, \qquad l, k=1,2 $$ and the constants
$\kappa_1^{(l)}$ are arbitrary. We shall define these constants
from the solvability conditions for $u^{(l)}_2$. Here and
throughout $l^*=1$ if $l=2$ and $l^*=2$ if $l=1.$

Thus,
\begin{equation*}\label{5000}
\alpha_{11}^{(l)}=\widetilde\alpha_{11}^{(l)}+\kappa_1^{(l)}
\alpha_{01}^{(l^*)},
\end{equation*}
where $$\widetilde \alpha_{11}^{(l)}=\frac{\partial \widetilde
u_1^{(l)}}{\partial x_2}\bigg|_{x_2=0},\qquad
\frac{d^{2k+1}\widetilde
\alpha_{11}^{(l)}}{dx_1^{2k+1}}\bigg|_{x_1
=\pm\frac{1}{2}}=0,\qquad k=0,1,\dots$$

Let us define $\alpha^{(l)}_{20}(x_1)$ and $v_2(\xi; x_2)$.
Consider an auxiliary problem:
\begin{equation}\label{6!}
\left\{\begin{array}{l} \Delta_\xi \widetilde{X}=\frac{\partial
{X}}{\partial \xi_1}\quad\hbox{in} \ \Pi,\\
\widetilde{X}=0\quad\hbox{on} \ \partial\Pi.
\end{array}\right.
\end{equation}
It is proved in \cite{AmCheGa} that Problem {\rm (\ref{6!})} has
a solution with the asymptotics
\begin{equation}\label{7!!}
\widetilde{X}(\xi)=0 \quad\hbox{as}\ \xi_2\to+\infty,
\end{equation}
up to exponentially small terms. Note that, due to the evenness
of the functions $F$, the solution $\widetilde{X}$ of Problem
(\ref{6!}) is odd in $\xi_1$, and thus has a $1$-periodic
extension in $\xi_1$ for which we keep the same notation
$\widetilde{X}$.

Then it is easy to see that, due to (\ref{6'}), (\ref{7'}),
(\ref{700}), (\ref{6!}) and (\ref{7!!}), the function
\begin{equation}\label{14}
v_2^{(l)}(\xi;x_1)=\alpha_{11}^{(l)}(x_1) X(\xi)
-2(\alpha_{01}^{(l)})^\prime(x_1)\widetilde{X}(\xi)
\end{equation}
is the $1$-periodic solution to Problem~(\ref{61}), which has the
asymptotics
\begin{equation*}\label{15}
v_2^{(l)}(\xi;x_1)=
\alpha_{11}^{(l)}(x_1)\left(\xi_2+C(F)\right)\quad\hbox{as
$\xi_2\to +\infty$}
\end{equation*}
up to exponentially small terms, and also which satisfies the
conditions~(\ref{800}), (\ref{601}) for
\begin{equation}\label{5001}
\alpha^{(l)}_{20}(x_1)=C(F)\left(\widetilde
\alpha_{11}^{(l)}+\kappa_1^{(l)}\alpha_{01}^{(l^*)}\right).
\end{equation}
Thus we defined $v_2^{(l)}$ and $\alpha^{(l)}_{20}$ up to
$\kappa_1^{(l)}$, which is unknown yet.

It is easy to verify that, due to (\ref{700}), (\ref{14}) and the
boundary condition from (\ref{6!}),
\begin{equation}\label{5+4}
\frac{\partial v_2^{(l)}}{\partial x_1}
\left(\xi;x_1\right)=0\quad\hbox{as $\xi=\frac{x}{\varepsilon}$
and $x_1=\pm\frac{1}{2}$}.
\end{equation}
Hence Problem (\ref{62'}) takes the form
\begin{equation}\label{62}
\left\{\begin{array}{l} -\Delta_\xi v_3^{(l)}= 2\frac{\partial^2
v_2^{(l)}}{\partial x_1\partial\xi_1}+
\frac{\partial^2 v_1^{(l)}}{\partial x_1^2}+\lambda_0 v_1^{(l)}\quad\hbox{in} \ \Pi,\\
v_3^{(l)}=0\quad\hbox{on} \ \Gamma,\\
\frac{\partial v_3^{(l)}}{\partial \xi_1}=0\quad\hbox{as} \
\xi_1=\pm\frac{1}{2},\ x_1=\pm \frac{1}{2}
.
\end{array}\right.
\end{equation}
and by (\ref{5+1}), (\ref{5+4}) and the boundary condition from
(\ref{61}), we have
\begin{equation}\label{5+5}
\frac{\partial v_2^{(l)}}{\partial \nu}
\left(\frac{x}{\varepsilon};x_1\right)=0\quad\hbox{on
$\Gamma_{2,\varepsilon}\cup\Gamma_{3,\varepsilon}$}.
\end{equation}
Let us define $\lambda_2^{(l)}$, $u_2^{(l)}$ and
$\kappa_1^{(l)}$. From (\ref{solva}), (\ref{5001}), (\ref{2p2})
and (\ref{6.26}) we deduce that the sufficient solvability
conditions of Problem (\ref{91}),~(\ref{93}) are
\begin{equation*}\label{6.22}
\lambda^{(l)}_2=-C(F)\int\limits_{-\frac{1}{2}}^{\frac{1}{2}}
\widetilde\alpha^{(l)}_{11}(x_1)\alpha^{(l)}_{01}(x_1)\ dx_1
\end{equation*}
and
\begin{equation*}\label{88}
\kappa_1^{(l)}=\frac{\displaystyle
\int\limits_{-\frac{1}{2}}^{\frac{1}{2}}
\widetilde\alpha^{(l)}_{11}(x_1)\alpha^{(l^*)}_{01}(x_1)\ dx_1}
{\displaystyle
\int\limits_{-\frac{1}{2}}^{\frac{1}{2}}\left(\left(\alpha^{(l)}_{01}\right)^2(x_1)-
\left(\alpha^{(l^*)}_{01}\right)^2(x_1)\right)\ dx_1}.
\end{equation*}
Thus, we defined constants $\lambda_2^{(l)}$ è $\kappa_1^{(l)}$
and in particular because of (\ref{ttr}) the function
$u_1^{(l)}(x)$. We choose $u^{(l)}_1$ in the form:
\begin{equation}\label{ttr1}
u^{(l)}_2=\widetilde u^{(l)}_2+\kappa_2^{(l)} u^{(l^*)}_0,
\end{equation}
where $$ \int\limits_{\Omega}\widetilde u^{(l)}_2(x)
u^{(k)}_0(x)\ dx=0, \qquad l, k=1,2 $$ and the constants
$\kappa_2^{(l)}$ are arbitrary. We shall define these constants
from the solvability conditions for $u^{(l)}_3$.

Thus,
\begin{equation*}\label{6000}
\alpha_{21}^{(l)}=\widetilde\alpha_{21}^{(l)}+\kappa_2^{(l)}
\alpha_{01}^{(l^*)},
\end{equation*}
where $$ \widetilde \alpha_{21}^{(l)}=\frac{\partial \widetilde
u_2^{(l)}}{\partial x_2}\bigg|_{x_2=0},\qquad
\frac{d^{2k+1}\widetilde
\alpha_{21}^{(l)}}{dx_1^{2k+1}}\bigg|_{x_1
=\pm\frac{1}{2}}=0,\qquad k=0,1,\dots. $$

Let us define $\alpha^{(l)}_{30}(x_1)$ and $v_3(\xi; x_2)$.
Consider auxiliary problems:
\begin{equation}\label{7!1}
\left\{\begin{array}{l}
\Delta_\xi\widetilde{\widetilde{X}}_{(I)}=\frac{\partial
{\widetilde{X}}}{\partial \xi_1}\quad\hbox{in} \ \Pi,\\
\widetilde{\widetilde{X}}_{(I)}=0\quad\hbox{on} \ \Gamma,\qquad
\frac{\partial \widetilde{\widetilde{X}}_{(I)}}{\partial
\xi_1}=0\quad\hbox{as} \ \xi_1=\pm\frac{1}{2}.
\end{array}\right.
\end{equation}
\begin{equation}\label{7!2}
\left\{\begin{array}{l}
\Delta_\xi\widetilde{\widetilde{X}}_{(II)}=X\quad\hbox{in} \ \Pi,\\
\widetilde{\widetilde{X}}_{(II)}=0\quad\hbox{on} \ \Gamma,\qquad
\frac{\partial \widetilde{\widetilde{X}}_{(II)}}{\partial
\xi_1}=0\quad\hbox{as} \ \xi_1=\pm\frac{1}{2}.
\end{array}\right.
\end{equation}
In Section~\ref{s4} we shall prove the following two statements.
\begin{proposition}\label{p1}
Problem {\rm (\ref{7!1})} has a solution with the asymptotics
\begin{equation}\label{88!!}
\widetilde{\widetilde{X}}_{(I)}(\xi)=C_{(I)}(F)  \quad\hbox{as}\
\xi_2\to+\infty,
\end{equation}
up to exponentially small terms, where $C_{(I)}(F)$ is a constant
depending on the function $F$.
\end{proposition}

\begin{proposition}\label{p2}
Problem {\rm(\ref{7!2})} has a solution with the asymptotics
\begin{equation}\label{89!!}
\widetilde{\widetilde{X}}_{(II)}(\xi)=\frac{1}{6}\xi_2^3+\frac{1}{2}C(F)\xi_2^2+C_{(II)}(F)
\quad\hbox{as}\ \xi_2\to+\infty,
\end{equation}
up to exponentially small terms, where $C_{(II)}(F)$ is a constant
depending on the function $F$.
\end{proposition}

Note that, due to the evenness of the function $F$ and the
evenness of the right-hand sides of the equations in~(\ref{7!1})
and~(\ref{7!2}) the solutions $\widetilde{\widetilde{X}}_{(I)}$
and $\widetilde{\widetilde{X}}_{(II)}$ of Problems (\ref{7!1})
and~(\ref{7!2}) respectively are even in $\xi_1$, and thus have
$1$-periodic extensions in $\xi_1$ for which we keep the same
notation $\widetilde{\widetilde{X}}_{(I)}$ and
$\widetilde{\widetilde{X}}_{(II)}$.

Then it is easy to see that, due to (\ref{6'}), (\ref{7'}),
(\ref{700}), (\ref{6!}), (\ref{7!!}), (\ref{7!1})--(\ref{89!!})
the function
\begin{equation}\label{140}
\begin{aligned}
v_3^{(l)}(\xi;x_1)&=
\alpha_{21}^{(l)}(x_1) X(\xi)
-2(\alpha_{11}^{(l)})^\prime(x_1)\widetilde{X}(\xi)+\\
&+4(\alpha_{01}^{(l)})^{\prime\prime}(x_1)\widetilde{\widetilde{X}}_{(I)}(\xi)
-(\alpha_{01}^{(l)})^{\prime\prime}(x_1)\widetilde{\widetilde{X}}_{(II)}(\xi)
-\lambda_0\alpha_{01}^{(l)}(x_1)\widetilde{\widetilde{X}}_{(II)}(\xi)
\end{aligned}
\end{equation}
is the $1$-periodic solution to Problem~(\ref{62}), which has the
asymptotics
\begin{equation}\label{150}
\begin{aligned}
v_3^{(l)}(\xi;x_1)&=\alpha_{21}^{(l)}(x_1)\left(\xi_2+C(F)\right)+4C_{(I)}(F)
(\alpha_{01}^{(l)})^{\prime\prime}(x_1)-\\
&-\left((\alpha_{01}^{(l)})^{\prime\prime}(x_1)+
\lambda_0\alpha_{01}^{(l)}(x_1)\right)
\left(\frac{1}{6}\xi_2^3+\frac{1}{2}C(F)\xi_2^2+C_{(II)}(F)\right)\\
&\hbox{as $\xi_2\to +\infty$,}
\end{aligned}
\end{equation}
up to exponentially small terms. Note that due to equations for
$u_k^{(l)}$ from (\ref{90}), (\ref{91}) we have: $$
-\frac{1}{6}\left((\alpha_{01}^{(l)})^{\prime\prime}(x_1)+
\lambda_0\alpha_{01}^{(l)}(x_1)\right)=\alpha_{03}^{(l)}(x_1) $$
and $$
-\frac{1}{2}C(F)\left((\alpha_{01}^{(l)})^{\prime\prime}(x_1)+
\lambda_0\alpha_{01}^{(l)}(x_1)\right)=\alpha_{12}^{(l)}(x_1). $$

Hence from (\ref{150}) we conclude that $v_3$ satisfies the
condition~(\ref{800}), (\ref{602}) for $$
\alpha^{(l)}_{30}(x_1)=C(F)\left(\widetilde
\alpha_{21}^{(l)}(x_1)+\kappa_2^{(l)}\alpha_{01}^{(l^*)}(x_1)\right)+4C_{(I)}(F)
(\alpha_{01}^{(l)})^{\prime\prime}(x_1)- $$
\begin{equation}\label{6002}
-C_{(II)}(F)\left((\alpha_{01}^{(l)})^{\prime\prime}(x_1)+
\lambda_0\alpha_{01}^{(l)}(x_1)\right).
\end{equation}

Thus we defined $v_3^{(l)}$ and $\alpha^{(l)}_{30}$ up to
$\kappa_2^{(l)}$, which is unknown yet.

It is easy to verify that, due to (\ref{700}), (\ref{140}) and the
boundary condition from (\ref{7!1}) and~(\ref{7!2}),
\begin{equation}\label{5+4+1}
\frac{\partial v_3^{(l)}}{\partial x_1}
\left(\xi;x_1\right)=0\quad\hbox{as $\xi=\frac{x}{\varepsilon}$
and $x_1=\pm\frac{1}{2}$},
\end{equation}
and hence, by (\ref{5+1}), (\ref{5+4+1}) and the boundary
condition from (\ref{62}), we have
\begin{equation}\label{5+5+1}
\frac{\partial v_3^{(l)}}{\partial \nu}
\left(\frac{x}{\varepsilon};x_1\right)=0\quad\hbox{on
$\Gamma_{2,\varepsilon}\cup\Gamma_{3,\varepsilon}$}.
\end{equation}
Let us define $\lambda_3^{(l)}$, $u_3^{(l)}$ and
$\kappa_2^{(l)}$. From (\ref{solva}), (\ref{6002}), (\ref{2p2})
and (\ref{6.26}) we deduce that the sufficient solvability
conditions of Problem (\ref{92}),~(\ref{93}) are $$
\lambda^{(l)}_3=\int\limits_{-\frac{1}{2}}^{\frac{1}{2}}\alpha^{(l)}_{30}(x_1)\alpha^{(l)}_{01}(x_1)\
dx_1=-C(F)\int\limits_{-\frac{1}{2}}^{\frac{1}{2}}\widetilde
\alpha_{21}^{(l)}(x_1)\alpha_{01}^{(l)}(x_1) dx_1+ $$
\begin{equation*}\label{6.220}
+\left(C_{(II)}(F)-4C_{(I)}(F)\right)
\int\limits_{-\frac{1}{2}}^{\frac{1}{2}}(\alpha^{(l)}_{01})^{\prime\prime}(x_1)
\alpha_{01}^{(l)}(x_1) dx_1\ dx_1+
\end{equation*}
$$ +\lambda_0 C_{(II)}(F)
\int\limits_{-\frac{1}{2}}^{\frac{1}{2}}\left(\alpha^{(l)}_{01}\right)^2(x_1)
dx_1 $$ and
$$
\lambda_1^{(l)}\kappa_2^{(l)}+\lambda_2^{(l)}\kappa_1^{(l)}=-C(F)\int\limits_{-\frac{1}{2}}^{\frac{1}{2}}\left(\widetilde
\alpha_{21}^{(l)}(x_1)\alpha_{01}^{(l^*)}(x_1)+\kappa_2^{(l)}
\left(\alpha^{(l^*)}_{01}\right)^2(x_1)\right) dx_1+
$$
\begin{equation*}\label{6.220'}
+\left(C_{(II)}(F)-4C_{(I)}(F)\right)
\int\limits_{-\frac{1}{2}}^{\frac{1}{2}}(\alpha^{(l)}_{01})^{\prime\prime}(x_1)
\alpha_{01}^{(l^*)}(x_1) \ dx_1.
\end{equation*}
From the last formula we deduce the expression for
$\kappa_2^{(l)}$. Thus, we defined constants $\lambda_3^{(l)}$,
$\kappa_2^{(l)}$ and in particular because of (\ref{ttr1}) the
function $u_2^{(l)}(x)$.
We fix the arbitrariness in choosing of $u^{(l)}_3$ by the
following $$ u^{(l)}_3=\widetilde u^{(l)}_3+\kappa_3^{(l)}
u^{(l^*)}_0, $$  where $$ \int\limits_{\Omega}\widetilde
u^{(l)}_3(x) u^{(k)}_0(x)\ dx=0, \qquad l,k=1,2. $$ Here
constants $\kappa_3^{(l)}$ are arbitrary and we shall define
them in a unique way from the solvability conditions for
$u^{(l)}_4$.

\begin{remark}\label{rem2}
In the same way we can construct the complete asymptotic
expansion of eigenelements of Problem (\ref{2p3}) in the form
(\ref{uu})--(\ref{ve}). Substituting (\ref{uu}) and~(\ref{-1})
in~(\ref{2p3}), we deduce the boundary value problems for
coefficients of~(\ref{uu}):
\begin{equation}\label{92i}
\left\{\begin{array}{l} -\Delta u_i^{(l)}=\lambda_0
u_i^{(l)}+\sum\limits_{k=1}^{i-1}\lambda_k^{(l)} u_{i-k}^{(l)}+
\lambda_i^{(l)}u_0^{(l)}\quad\hbox{in} \ \Omega,\\
\frac{\partial u_i^{(l)}}{\partial \nu}=0\quad\hbox{on} \
\Gamma_1\cup\Gamma_2\cup\Gamma_3,\\ u_i^{(l)}=
\alpha_{i0}^{(l)}\quad\hbox{on} \ \Gamma_0.
\end{array}\right.
\end{equation}
We construct the solution of Problem (\ref{92i}) in the form
\begin{equation}\label{ttri}
u^{(l)}_i=\widetilde u^{(l)}_i+\kappa_i^{(l)} u^{(l^*)}_0,
\end{equation}
where
\begin{equation}\label{oop}
\int\limits_{\Omega}\widetilde u^{(l)}_i(x) u^{(k)}_0(x)\ dx=0,
\qquad l,k=1,2.
\end{equation}
Rewriting the asymptotics of~(\ref{uu}) as $x_2\to0$ in the fast
variables $\xi=\frac{x}{\varepsilon}$ we obtain that
\begin{equation}\label{900i}
\sum\limits_{i=0}^{\infty}\varepsilon^i
u_i^{(l)}(x)=\sum\limits_{i=1}^{\infty}\varepsilon^i
V_i^{(l)}(\xi; x_1)\qquad\hbox{as}\ \ x_2=\varepsilon\xi_2\to0,
\end{equation}
where
\begin{equation*}
V_i^{(l)}(\xi; x_1)=\widetilde{V}_i^{(l)}(\xi; x_1)+
\left(\widetilde{\alpha}^{(l)}_{i-1,1}(x_1)+\kappa_{i-1}^{(l)}
\alpha^{(l^*)}_{0,1}(x_1)\right)\xi_2
+\alpha^{(l)}_{i0}(x_1)
\end{equation*}
and $\widetilde{V}_i^{(l)}(\xi; x_1)$ is independent of $u_j$
for $j\geq i-1$, $\widetilde{\alpha}^{(l)}_{i-1,1}(x_1)$ depends
only on $\widetilde{u}_{i-1}$. Substituting~(\ref{-1}),
(\ref{ve}) in~(\ref{2p3}) keeping in mind~(\ref{900i}), we get
the equations and the boundary conditions for the terms
$v_i^{(l)}$ of (\ref{ve}):
\begin{equation}\label{62i}
\left\{\begin{array}{l} -\Delta_\xi v_i^{(l)}= 2\frac{\partial^2
v_{i-1}^{(l)}}{\partial x_1\partial\xi_1}+ \frac{\partial^2
v_{i-2}^{(l)}}{\partial x_1^2}+\sum\limits_{k=0}^{i-3}\lambda_k
v_{i-2-k}^{(l)}\quad\hbox{in} \ \Pi,\\ v_i^{(l)}=0\quad\hbox{on}
\ \Gamma,\\ \frac{\partial v_i^{(l)}}{\partial
\xi_1}=-\frac{\partial v_{i-1}^{(l)}}{\partial
x_1}\quad\hbox{as} \ \xi_1=\pm\frac{1}{2},\ x_1=\pm
\frac{1}{2},\\ v_i^{(l)}\sim V_i^{(l)} \quad\hbox{as}\
\xi_2\to+\infty
.
\end{array}\right.
\end{equation}
Before the $i$-th step we have defined $\lambda_k^{(l)}$,
$v_k^{(l)}$, $k\leq i-1$, $u_{j}^{(l)}$, $j\leq i-2$,
$\widetilde{u}_{i-1}^{(l)}$ (and consequently
$\widetilde{V}_i^{(l)}$ and
$\widetilde{\alpha}^{(l)}_{i-1,1}(x_1)$).

On the $i$-th step solving Problem (\ref{62i}), we find the
boundary condition $\alpha_{i0}^{(l)}$ in the form
\begin{equation}\label{ppo}
\alpha_{i0}^{(l)}=\widetilde{\alpha}_{i0}^{(l)}+
\kappa_i^{(l)}C(F)\alpha_{01}^{(l^*)},
\end{equation}
where $\widetilde{\alpha}_{i0}^{(l)}$ depends only on
$\widetilde{V}_i^{(l)}$ and
$\widetilde{\alpha}^{(l)}_{i-1,1}(x_1)$. Then from the
solvability condition for Problem (\ref{92i}) with
$\alpha_{i0}^{(l)}$ defined in (\ref{ppo}),  we derive
$\lambda_i^{(l)}$ and $\kappa_{i-1}^{(l)}$ (hence exactly define
$u_{i-1}^{(l)}$). We choose $u^{(l)}_i$ in the form
(\ref{ttri}), (\ref{oop}), where the constants $\kappa_i^{(l)}$
are arbitrary. We shall define these constants in the next step.
Thus we completed the $i$-th step.

Acting in the same way, we construct complete formal asymptotic
expansions (\ref{uu})--(\ref{ve}) of eigenelements.

On the justification of these asymptotics see Remark \ref{rem3}.
\end{remark}

\section{Verification of the asymptotics.}\label{s3}

Denote
\begin{equation}\label{2.2}
\widetilde \lambda_\varepsilon^{(l)}=\lambda_0+\varepsilon
\lambda_1^{(l)}+\varepsilon^2\lambda_2^{(l)}+
\varepsilon^3\lambda_3^{(l)},
\end{equation}
\begin{equation}
\begin{aligned}
\widetilde u_\varepsilon^{(l)}(x) =&
\left(u_0^{(l)}(x)+\sum\limits_{i=1}^3\varepsilon^i
u_i^{(l)}(x)\right) \chi\left(\frac{x_2}{\varepsilon^\beta}\right)
\\
&+\left(\sum\limits_{i=1}^3\varepsilon^i
v_i^{(l)}\left(\frac{x}{\varepsilon};x_1\right)\right)
\left(1-\chi\left(\frac{x_2} {\varepsilon^\beta}\right)\right),
\end{aligned}\label{2.1}
\end{equation}
where $\chi(s)$ is a smooth cut-off function, equals to zero as
$s<1$ and equals to one as $s>2$, and $\beta$ is a fixed number
($0<\beta<1$). Obviously, $\widetilde u_\varepsilon^{(l)}\in
C^\infty(\overline{\Omega})$

\begin{lemma}\label{lm3.1}
The function $\widetilde u_\varepsilon^{(l)}$ is the
solution of Problem
\begin{align}\label{4.01'}
\left\{\begin{array}{l} -\Delta \widetilde
u_\varepsilon^{(l)}=\widetilde{\lambda}_\varepsilon \widetilde
u_\varepsilon^{(l)}+f^{(l)}_\varepsilon\quad\hbox{in} \
\Omega^\varepsilon,\\ \widetilde
u_\varepsilon^{(l)}=0\quad\hbox{on} \ \Gamma_\varepsilon,\qquad
\frac{\partial \widetilde u_\varepsilon^{(l)}}{\partial
\nu}=0\quad\hbox{on} \
\Gamma_1\cup\Gamma_{2,\varepsilon}\cup\Gamma_{3,\varepsilon},
\end{array}\right.
\end{align}
where
\begin{equation}\label{IF}
\left\|f_\varepsilon^{(l)}\right\|_{L_2(\Omega^\varepsilon)}=O(\varepsilon^{5
\beta/2}).
\end{equation}
\end{lemma}
\begin{proof}
Due to boundary conditions from (\ref{2}), (\ref{90}),~(\ref{91})
and~(\ref{92}) for $u_0^{(l)}$, $u_1^{(l)}$, $u_2^{(l)}$,
$u_3^{(l)}$ and the boundary conditions (\ref{5+3}), (\ref{5+5}),
(\ref{5+5+1}) for $v_1^{(l)}$, $v_2^{(l)}$, $v_3^{(l)}$, we get
that the function $\widetilde{u}_\varepsilon^{(l)}$ satisfies the
boundary conditions of Problem~(\ref{4.01'}).

From the other hand, due to  Formula (\ref{5}), and equations
(\ref{2}), (\ref{90}), (\ref{91}), (\ref{92}),
(\ref{60}),~(\ref{61}) and (\ref{62}) for $u_0^{(l)}$,
$u_1^{(l)}$, $u_2^{(l)}$, $u_3^{(l)}$, $v_1^{(l)}$, $v_2^{(l)}$
and $v_3^{(l)}$, the function $\widetilde{u}_\varepsilon^{(l)}$
satisfies the equation of Problem~(\ref{4.01'}) where
\begin{equation*}
-f_\varepsilon^{(l)}(x)= I_1^{(l)}(x;\varepsilon)+
I_2^{(l)}(x;\varepsilon)+
 I_3^{(l)}(x;\varepsilon)
+I_4^{(l)}(x;\varepsilon),
\end{equation*}
\begin{equation*}
\begin{aligned}
I_1^{(l)}&=\varepsilon^4\chi\left(\frac{x_2}
{\varepsilon^\beta}\right)\left(\lambda_1^{(l)}
u_3^{(l)}+\lambda_2^{(l)} u_2^{(l)}+\varepsilon\lambda_2^{(l)}
u_3^{(l)}+\lambda_3^{(l)} u_1^{(l)}+\varepsilon\lambda_3^{(l)}
u_2^{(l)}+\varepsilon^2\lambda_3^{(l)} u_3^{(l)}\right), \\
I_2^{(l)}&=\varepsilon^2\left(1-\chi\left(\frac{x_2}
{\varepsilon^\beta}\right)\right)\Bigg(\left(\lambda_0
v_2^{(l)}+\varepsilon\lambda_0 v_3^{(l)}+\lambda_1
v_1^{(l)}+\varepsilon\lambda_1^{(l)}
v_2^{(l)}+\varepsilon^2\lambda_1
v_3^{(l)}+\varepsilon\lambda_2^{(l)}
v_1^{(l)}\right)\\
&+\varepsilon^2\left(\varepsilon^2\lambda_2^{(l)}
v_2^{(l)}+\varepsilon^2\lambda_3^{(l)}
v_1^{(l)}+\varepsilon^3\lambda_3^{(l)}
v_2^{(l)}+\varepsilon^4\lambda_3^{(l)} v_3^{(l)}\right)\\&
+\varepsilon^2\left(\frac{\partial^2 v_2^{(l)}}{\partial
x_1^2}+\varepsilon\frac{\partial^2 v_3^{(l)}}{\partial x_1^2}
+2\frac{\partial^2 v_3^{(l)}}{\partial
x_1\partial\xi_1}\right)\Bigg),\\
I_3^{(l)}&=2\varepsilon^{-\beta} \chi'\left(\frac{x_2}
{\varepsilon^\beta}\right)\left(\sum\limits_{i=0}^{3}\varepsilon^i
\frac{\partial u_i^{(l)}}{\partial
x_2}-\sum\limits_{i=1}^{3}\varepsilon^{i-1}\frac{\partial
v_i^{(l)}}{\partial \xi_2}\right),\\
I_4^{(l)}&=\varepsilon^{-2\beta} \chi''\left(\frac{x_2}
{\varepsilon^\beta}\right) \left( u_0^{(l)}+\varepsilon
u_1^{(l)}+\varepsilon^2 u_2^{(l)}+\varepsilon^3
u_3^{(l)}-(\varepsilon v_1^{(l)}+\varepsilon^2
v_2^{(l)}+\varepsilon^3 v_3^{(l)})\right).
\end{aligned}
\end{equation*}
Since the functions $u_i^{(l)}$ are smooth, then it is obvious
that
\begin{equation}\label{I1}
\left\|I_1^{(l)}\right\|_{L_2(\Omega^\varepsilon)}=O(\varepsilon^4).
\end{equation}
Due to (\ref{1!}), (\ref{7}) and (\ref{14}), we obtain that
\begin{equation}\label{I2}
\left\| I_2^{(l)}\right\|_{L_2(\Omega^\varepsilon)}=
O\left(\varepsilon^{\frac{5\beta}{2}}\right).
\end{equation}
Bearing in mind the matching conditions  (\ref{900}), (\ref{800})
and that the derivatives of $\chi\left(\frac{x_2}
{\varepsilon^\beta}\right)$ are not equal to zero only in the
strip $\varepsilon^\beta<x_2<2\varepsilon^\beta$, it is easy to
see that
\begin{equation}\label{I}
\left\|I_3^{(l)}\right\|_{L_2(\Omega^\varepsilon)}+
\left\|I_4^{(l)}\right\|_{L_2(\Omega^\varepsilon)}=
O(\varepsilon^{\frac{5\beta}{2}}).
\end{equation}
From (\ref{I1})--(\ref{I}) it follows (\ref{IF}).
\end{proof}

The following statement is proved in \cite{AmCheGa}.
\begin{lemma}\label{lm3.3} Assume that
the multiplicity of the eigenvalue $\lambda_0$ of Problem
(\ref{2}) is equal to $p$. Then for any $\lambda$ close to
$\lambda_0$
\begin{enumerate}
\item[(i)]  the solution
$U_\varepsilon$ to Problem
\begin{equation}\label{4.01}
\left\{\begin{array}{l} -\Delta U_\varepsilon=\lambda
U_\varepsilon+F_\varepsilon\quad\hbox{in} \ \Omega^\varepsilon,\\
U_\varepsilon=0\quad\hbox{on} \ \Gamma_\varepsilon,\qquad
\frac{\partial U_\varepsilon}{\partial \nu}=0\quad\hbox{on} \
\Gamma_1\cup\Gamma_{2,\varepsilon}\cup\Gamma_{3,\varepsilon}
\end{array}\right.
\end{equation}
admits the estimate
\begin{equation}\label{4.14}
\|U_\varepsilon\|_{H^1(\Omega^\varepsilon)}\leq
\mathcal{C}\frac{\|F_\varepsilon\|_{L_2(\Omega^\varepsilon)}}
{\prod\limits_{l=1}^p |{\lambda^{(l)}_\varepsilon}-\lambda|}
\end{equation}
where $\lambda_\varepsilon^{(1)},\dots,\lambda_\varepsilon^{(p)}$
are the eigenvalues of Problem {\rm(\ref{1})}, which converge to
$\lambda_0$;
\item[(ii)] if a solution $U_\varepsilon$ to Problem {\rm
(\ref{4.01})} is orthogonal in $L_2(\Omega^\varepsilon)$ to the
eigenfunction $u^{(k)}_\varepsilon$ of Problem~{\rm(\ref{1})}
corresponding to $\lambda_\varepsilon^{(k)}$, then it satisfies
the estimate
\begin{equation}\label{e2}
\|U_\varepsilon\|_{H^1(\Omega^\varepsilon)}\leq
\mathcal{C}\frac{\|F_\varepsilon\|_{L_2(\Omega^\varepsilon)}}
{\prod\limits_{l=1;\,l\not=k}^p
|{\lambda^{(l)}_\varepsilon}-\lambda|}.
\end{equation}
\end{enumerate}
\end{lemma}
For our case $p=2$ from this lemma we deduce the statement.
\begin{corollary}\label{cor3.3} For any $\lambda$ close to $\lambda_0$
\begin{enumerate}
\item[(i)]  the solution
$U_\varepsilon$ to Problem {\rm(\ref{4.01})} admits the estimate
\begin{equation}\label{4.14'}
\|U_\varepsilon\|_{H^1(\Omega^\varepsilon)}\leq
\mathcal{C}\frac{\|F_\varepsilon\|_{L_2(\Omega^\varepsilon)}}
{\left(\max\limits_{l}
|\lambda^{(l)}_\varepsilon-\lambda|\right)^2},
\end{equation}
where $\lambda_\varepsilon^{(1)},\lambda_\varepsilon^{(2)}$ are
the eigenvalues of Problem {\rm(\ref{1})}, which converge to
$\lambda_0$;
\item[(ii)] if a solution $U_\varepsilon$ to Problem {\rm
(\ref{4.01})} is orthogonal in $L_2(\Omega^\varepsilon)$ to the
eigenfunction $u^{(l)}_\varepsilon$ of Problem~{\rm(\ref{1})}
corresponding to $\lambda_\varepsilon^{(l)}$, then it satisfies
the estimate
\begin{equation}\label{e2'}
\|U_\varepsilon\|_{H^1(\Omega^\varepsilon)}\leq
\mathcal{C}\frac{\|F_\varepsilon\|_{L_2(\Omega^\varepsilon)}} {
|{\lambda^{(k)}_\varepsilon}-\lambda|}\qquad\hbox{for $k\not=l$}.
\end{equation}
\end{enumerate}
\end{corollary}
\begin{theorem}\label{th3.1}
Assume that the multiplicity of $\lambda_0$ of Problem
{\rm(\ref{2})} equals two, the associated eigenfunctions
$u_0^{(l)}$ ($l=1, 2$) satisfy the conditions
{\rm(\ref{2p2})--(\ref{nerav})}. Then eigenvalues
$\lambda_\varepsilon^{(l)}$  of Problem {\rm(\ref{1})}, converging
to $\lambda_0$ as $\varepsilon \to 0$, and the associated
eigenfunctions $u_\varepsilon^{(l)}$ orthonormalized in
$L_2(\Omega_\varepsilon)$ have the following asymptotics:
\begin{align}\label{est1'}
\lambda_\varepsilon^{(l)}=&\widetilde{\lambda}^{(l)}_\varepsilon
+o(\varepsilon^{\frac{5\beta}{4}})\qquad\hbox{for any
$\beta<1$},\\
&\|u^{(l)}_\varepsilon-\widetilde{u}^{(l)}_\varepsilon\|_{H^1(\Omega^\varepsilon)}
=o(1).\label{Est2'}
\end{align}
\end{theorem}
\begin{proof}
Since $u_j^{(l)}\in C^\infty(\overline{\Omega})$, then due to
(\ref{2.1}) and (\ref{600})--(\ref{800}) we derive
\begin{equation}\label{IF1}
\left\|\widetilde{u}_\varepsilon^{(l)}\right\|_{L_2(\Omega^\varepsilon)}=1+o(1)\qquad\hbox{as
$\varepsilon\to0$}.
\end{equation}

Applying item {\it (i)} of Corollary~\ref{cor3.3} for
$\lambda=\widetilde{\lambda}^{(l)}_\varepsilon$,
$F_\varepsilon=f^{(l)}_\varepsilon$ and
$U_\varepsilon=\widetilde{u}^{(l)}_\varepsilon$, due to
Lemma~\ref{lm3.1} we obtain (\ref{est1'}).

Denote
\begin{equation}
\widehat{u}_\varepsilon^{(l)}=\widetilde{u}_\varepsilon^{(l)}-
\left(\widetilde{u}_\varepsilon^{(l)},u_\varepsilon^{(l)}\right)
_{L_2(\Omega^\varepsilon)}{u}_\varepsilon^{(l)}. \label{8.26}
\end{equation}
By the definition
\begin{equation}
\left(\widehat{u}_\varepsilon^{(l)},
{u}_\varepsilon^{(l)}\right)_{L_2(\Omega^\varepsilon)}=0,
\label{8.27}
\end{equation}
and the function $U_\varepsilon=\widehat{u}_\varepsilon^{(l)}$ is
a solution of Problem (\ref{4.01}) for
$\lambda=\widetilde{\lambda}^{(l)}_\varepsilon$ and
\begin{equation}
F_\varepsilon={f}^{(l)}_{\varepsilon}+\left(
\widetilde{\lambda}^{(l)}_
\varepsilon-\lambda^{(l)}_\varepsilon\right)
\widetilde{u}^{(l)}_{\varepsilon}. \label{8.28}
\end{equation}
From (\ref{8.28}), (\ref{est1'}), (\ref{IF1}) and (\ref{IF}) it
follows that
\begin{equation}\label{IF'}
\left\|F_\varepsilon\right\|_{L_2(\Omega^\varepsilon)}=
O(\varepsilon^{\frac{5\beta}{4}}).
\end{equation}
Since due to (\ref{est1}) and (\ref{nerav}) we have
$\left|\lambda^{(1)}_\varepsilon-
\lambda^{(2)}_\varepsilon\right|> c\varepsilon$, where $c>0$, then
by (\ref{IF'}) and item {\it (ii)} of Corollary~\ref{cor3.3} it
follows that
\begin{equation}
\left\| \widehat{u}_{\varepsilon}^{(l)}\right\|
_{H^1\left(\Omega^\varepsilon\right)}=o\left(1\right).
\label{8.29}
\end{equation}
Due to (\ref{8.26}), (\ref{8.29}) and (\ref{IF1}) we deduce
(\ref{Est2'}).
\end{proof}

{\it Proof of Theorem~\ref{t2}}. Since $u_j^{(l)}\in
C^\infty(\overline{\Omega})$, then due to (\ref{2.1}) and
(\ref{600})--(\ref{800}) we derive
\begin{equation}\label{Est2''}
\|\widetilde{u}^{(l)}_\varepsilon-u^{(l)}_0\|_{H^1(\Omega)} +
\|\widetilde{u}^{(l)}_\varepsilon\|_{H^1(\Omega^\varepsilon
\backslash\overline{\Omega})} =o(1).
\end{equation}
Then Theorem~\ref{t2} follows from Theorem~\ref{th3.1}.

\begin{remark}\label{rem3}
In an analogues way we can justify the complete asymptotic
expansion of the eigenelements of Problem (\ref{2p3})
constructed in Remark~\ref{rem2}. Denoting
\begin{equation}\label{2.2"}
\begin{aligned}
\widetilde
\lambda_\varepsilon^{(l)}=&\lambda_0+\sum\limits_{i=1}^{n}\varepsilon^i
\lambda_i^{(l)},\\
\widetilde u_\varepsilon^{(l)}(x) =&
\left(u_0^{(l)}(x)+\sum\limits_{i=1}^n\varepsilon^i
u_i^{(l)}(x)\right) \chi\left(\frac{x_2}{\varepsilon^\beta}\right)
\\
&+\left(\sum\limits_{i=1}^n\varepsilon^i
v_i^{(l)}\left(\frac{x}{\varepsilon};x_1\right)\right)
\left(1-\chi\left(\frac{x_2} {\varepsilon^\beta}\right)\right)
\end{aligned}
\end{equation}
and repeating the proof of Lemma~\ref{lm3.1}, we conclude that
Lemma~\ref{lm3.1} holds true for
\begin{equation}\label{IFI}
\left\|f_\varepsilon^{(l)}\right\|_{L_2(\Omega^\varepsilon)}=
O(\varepsilon^{N})\qquad\hbox{where
$N\underset{n\to\infty}\longrightarrow\infty$.}
\end{equation}
Obviously that (\ref{Est2''}) also holds true for $\widetilde
u_\varepsilon^{(l)}(x)$ defined in (\ref{2.2"}).  Then taking
into account (\ref{IFI}) and (\ref{Est2''}) and repeating the
proof of Theorem~\ref{th3.1}, we obtain that
\begin{align}\nonumber
\lambda_\varepsilon^{(l)}=&\widetilde{\lambda}^{(l)}_\varepsilon
+O(\varepsilon^{\frac{N}{2}})\qquad\hbox{where
$N\underset{n\to\infty}\longrightarrow\infty$.}
\end{align}
\end{remark}

\section{Appendix}\label{s4}

The proofs of Propositions~\ref{p1} and \ref{p2} are similar to
that of L.Tartar \cite{Tartar} (Lemma V.9) for a problem in a
semi-infinite strip whit a flat bottom (see also \cite{Lions}).
\begin{lemma}\label{l2}
Let $E$ and $E_0$ be Hilbert spaces, let $a$ be a continuous
bilinear form from $E\times E_0$, and let $M$ be a continuous
linear mapping from $E$ onto $E_0$. Assume that there exists
$\gamma>0$ such that $$a(u, Mu) \ge \gamma \Vert u \Vert^2
\quad\hbox{ for every } u\in E,$$ where $\Vert \cdot \Vert$
denotes the norm in $E$. Then, for every continuous linear form
$L$ into $E_0$, there exists a unique $u\in E$ satisfying $$a(u,
v)=L(v) \quad\hbox{ for every } v\in E_0.$$
\end{lemma}

{\it Proof of Proposition} \ref{p1}. In view of \cite{AmCheGa}
(proof of Proposition~1) there exists a positive constant
$\varsigma>0$ such that, $\forall \delta>0$, $\forall \alpha \in
\mathbb{N}^2$,
\begin{equation}\label{7!!!!}
\left\vert{\partial^\alpha \left(\frac{\partial
\widetilde{X}}{\partial \xi_1}\right)}
(\xi_1,\xi_2)\right\vert\leq C_{\delta,\alpha}\,
e^{-\varsigma\xi_2}
\end{equation}
for any $(\xi_1,\xi_2)\in\Pi$ with $\xi_2\ge \delta$, where
$C_{\delta,\alpha}$ is a constant depending only $\delta$ and
$\alpha$.
Let us introduce the Hilbert spaces $$E_\varsigma=\{v \ : \
e^{\varsigma\xi_2} v \in L_2(\Pi), \, e^{\varsigma\xi_2}
\frac{\partial v}{\partial \xi_j}\in L_2(\Pi) \quad\hbox{for}\
j=1, 2, \ v=0 \quad\hbox{on}\ \Gamma\},$$ $$E_\varsigma^0=\{v \
: \ v \in E_\varsigma, \, e^{\varsigma\xi_2} v\in L_2(\Pi)\},$$
equipped, respectively, with the scalar products (and associated
norms) $$(v,w)_\varsigma=\int\limits_\Pi e^{2 \varsigma\xi_2}
\nabla v \cdot \nabla w \ d\xi,$$
$$(v,w)_\varsigma^0=\int\limits_\Pi e^{2 \varsigma\xi_2}  v \,w
\ d\xi + \int\limits_\Pi e^{2 \varsigma\xi_2} \nabla v \cdot
\nabla w \ d\xi.$$ We consider the bilinear form $a_\varsigma$,
continuous on $E_\varsigma \times E_\varsigma^0$,
$$a_\varsigma(v,w)=\int\limits_\Pi \nabla v \cdot \nabla (e^{2
\varsigma\xi_2} w) \ d\xi \quad\hbox{for}\ v\in E_\varsigma, \,
w \in E_\varsigma^0,$$ and the linear form $L_\varsigma$,
continuous on $E_\varsigma^0$, $$L_\varsigma(v)=-\int\limits_\Pi
\frac{\partial \widetilde{X}} {\partial \xi_1} e^{2
\varsigma\xi_2} v \ d\xi \quad\hbox{for}\ v \in E_\varsigma^0.$$
Note here that, due to (\ref{7!!!!}), $\frac{\partial
\widetilde{X}} {\partial \xi_1} \in E_\varsigma^0$ and then
$L_\varsigma$ is well-defined.
We extend any $v \in V_\varsigma$ by 0 on $\{\xi\in
\mathbb{R}^2\
:\ -\frac{1}{2}<\xi_1<\frac{1}{2}, \ \xi_2<F(\xi_1)\}$, and we use
the same notation $v$ for the extension. For $v \in
E_\varsigma$, we denote
$$\overline{v}(\xi_2)=\int\limits_{-\frac{1}{2}}^{\frac{1}{2}}
v(\xi_1,\xi_2) \ d\xi_1, \quad\hbox{ for } \xi_2>0.$$ We have
the Friedrichs--Poincar\'e inequality
\begin{equation}\label{pin}
\int\limits_{-\frac{1}{2}}^{\frac{1}{2}} \left\vert
v(\xi_1,\xi_2)\right\vert^2 \ d\xi_1 \le \frac{1}{2}
\int\limits_{-\frac{1}{2}}^{\frac{1}{2}} \left\vert
\frac{\partial v }{\partial \xi_1} (\xi_1,\xi_2)\right\vert^2 \
d\xi_1 \quad\hbox{for} \ \xi_2<0,
\end{equation} and the Poincar\'e--Wirtinger
inequality
\begin{equation}\label{pwin}
\int\limits_{-\frac{1}{2}}^{\frac{1}{2}} \left\vert
v(\xi_1,\xi_2)-\overline{v}(\xi_2)\right\vert^2 \ d\xi_1 \le
\frac{1}{2} \int\limits_{-\frac{1}{2}}^{\frac{1}{2}} \left\vert
\frac{\partial v }{\partial \xi_1} (\xi_1,\xi_2)\right\vert^2 \
d\xi_1 \quad\hbox{for} \ \xi_2>0.
\end{equation}
Let $M: E_\varsigma \rightarrow E_\varsigma^0$ defined by
\begin{equation*}
M v=\left\{\begin{array}{l} v-\widetilde{v}(\xi_2) \quad\hbox{
in } \Pi_+,\\ v \quad\hbox{ in } \Pi_-,
\end{array}\right.
\end{equation*}
where $$\Pi^+=\{\xi\in \Pi : \xi=(\xi_1,\xi_2), \ \xi_2>0\},
\quad \Pi^-=\{\xi\in \Pi : \xi=(\xi_1,\xi_2), \ \xi_2<0\},$$ and
$\widetilde{v}$ is the solution of the differential equation
\begin{equation}
\left\{\begin{array}{l} \displaystyle{\frac{d \widetilde{v}}{d
\xi_2}} + 2 \varsigma \,\widetilde{v}=2 \varsigma \,\overline{v}
\quad\hbox{ for } \xi_2>0,\\ \widetilde{v}(0)=0.
\end{array}\right.\label{ed}
\end{equation}
We easily verify that
$e^{\varsigma\xi_2}(\widetilde{v}-\overline{v})$ solves the
differential equation
\begin{equation}
\displaystyle{\frac{d}{d
\xi_2}\left(e^{\varsigma\xi_2}(\widetilde{v}-\overline{v})\right)}
+ \varsigma \,e^{\varsigma\xi_2}
\,\left(\widetilde{v}-\overline{v}\right)= -
e^{\varsigma\xi_2}\,\frac{d \overline{v}}{d \xi_2} \quad\hbox{
for } \xi_2>0.\label{ed1}
\end{equation}
Since $e^{\varsigma\xi_2} \frac{\partial v}{\partial \xi_2} \in
L_2(\Pi)$ it follows that $e^{\varsigma\xi_2} \frac{d
\overline{v}}{d \xi_2} \in L_2(0,\infty)$. Multiplying
\eqref{ed1} by $e^{\varsigma\xi_2}(\widetilde{v}-\overline{v})$
and integrating on $(0,\infty)$ and using the trace theorem, it
follows that
\begin{equation}
\Vert e^{\varsigma\xi_2}(\widetilde{v}-\overline{v})
\Vert_{L_2(0,\infty)} \le C_\varsigma \Vert v
\Vert_{E_\varsigma} \quad \forall v \in E_\varsigma,
\label{ed11}
\end{equation}
where $C_\varsigma$ is a constant (depending on $\varsigma$). By
virtue of \eqref{pwin}, we have
\begin{equation}
\Vert e^{\varsigma\xi_2}(v-\overline{v}) \Vert_{L_2(\Pi^+)}^2
\le \frac{1}{2} \Vert v \Vert_{E_\varsigma}^2 \quad \forall v
\in E_\varsigma. \label{ed12}
\end{equation}
Using \eqref{ed11} and \eqref{ed12} we thus obtain that $M v \in
E_\varsigma^0$ and that the mapping $M$ is continuous from
$E_\varsigma$ to $E_\varsigma^0$. Let us verify that $M$ maps
$E_\varsigma$ onto $E_\varsigma^0$. Given $v$ in
$E_\varsigma^0$, we want to find $u\in E_\varsigma$ such that $M
u =v$. We seek $u$ in the form
\begin{equation*}
u=\left\{\begin{array}{l} v+h(\xi_2) \quad\hbox{ in } \Pi_+,\\ v
\quad\hbox{ in } \Pi_-.
\end{array}\right.
\end{equation*}
This implies $\overline{u}=\overline{v}+h \;$ (for $\xi_2>0)$,
then $h=\widetilde{u}$ and therefore $h$ is solution of the
differential equation
\begin{equation*}
\left\{\begin{array}{l} \displaystyle{\frac{d h}{d \xi_2}} + 2
\varsigma \,h=2 \varsigma (\overline{v}+h) \quad\hbox{ for }
\xi_2>0,\\ h (0)=0,
\end{array}\right.
\end{equation*}
that is $\displaystyle{\frac{d h}{d \xi_2}} =2 \varsigma
\,\overline{v}$. Since $v\in E_\varsigma^0$ it follows that
$e^{\varsigma\xi_2} \ \overline{v} \in L_2(0,\infty)$, and then
$e^{\varsigma\xi_2} \ {\frac{d h}{d \xi_2}} \in L_2(0,\infty)$.
Thus $u\in E_\varsigma$ and $M$ is onto.
Let us now prove that there exists a number $\gamma>0$ such that
$$a_\varsigma(v, Mv) \ge \gamma \Vert v \Vert^2 \quad\hbox{ for
every }v \in E_\varsigma.$$ We have
\begin{equation*}
\begin{array}{l} \displaystyle{a_\varsigma(v, M v) =\int\limits_\Pi \nabla v
\cdot \nabla (e^{2 \varsigma\xi_2} M v) \ d\xi}\\ =
\displaystyle \int\limits_{\Pi^+} \nabla v \cdot \nabla (e^{2
\varsigma\xi_2} (v-\widetilde{v})) \ d\xi + \int\limits_{\Pi^-}
\nabla v \cdot \nabla (e^{2 \varsigma\xi_2} v) \ d\xi\\
\displaystyle{= \int\limits_{\Pi^+} e^{2 \varsigma\xi_2} \vert
\nabla v\vert^2 \ d\xi +\int\limits_{\Pi^+} e^{2 \varsigma\xi_2}
\frac{\partial v }{\partial \xi_2} \left(2\,\varsigma
(v-\widetilde{v})-\frac{d \widetilde{v}}{d \xi_2}\right) \ d\xi}
\\ \displaystyle{+\int\limits_{\Pi^-} e^{2 \varsigma\xi_2} \vert
\nabla v\vert^2 \ d\xi +2 \varsigma \int\limits_{\Pi^-} e^{2
\varsigma\xi_2} \frac{\partial v }{\partial \xi_2} v \ d\xi}.
\end{array}
\end{equation*}
Using (\ref{ed}), it follows that
\begin{equation*}
\begin{array}{l}
a_\varsigma(v, M v) \displaystyle{=\int\limits_{\Pi} e^{2
\varsigma\xi_2} \vert \nabla v\vert^2 \ d\xi +2 \varsigma
\int\limits_{\Pi^+} e^{2 \varsigma\xi_2} \frac{\partial v
}{\partial \xi_2} (v-\overline{v}) \ d\xi}\\ \displaystyle{+ 2
\varsigma \int\limits_{\Pi^-} e^{2 \varsigma\xi_2}
\frac{\partial v }{\partial \xi_2} v \ d\xi.}
\end{array}
\end{equation*}
Applying the Young inequality and (\ref{pin}) and (\ref{pwin}),
it follows that
\begin{equation*}
\begin{array}{l}
a_\varsigma(v, M v) \displaystyle{\ge\int\limits_{\Pi} e^{2
\varsigma\xi_2} \vert \nabla v\vert^2 \ d\xi - \varsigma
\int\limits_{\Pi} e^{2 \varsigma\xi_2} \left\vert \frac{\partial
v }{\partial \xi_2}\right\vert^2 \ d\xi}\\
\displaystyle{-\frac{\varsigma}{2} \int\limits_{\Pi} e^{2
\varsigma\xi_2} \left\vert\frac{\partial v }{\partial
\xi_1}\right\vert^2 \ d\xi}  \\ \displaystyle{\ge \left(1-
\frac{3 k}{2}\right) \int\limits_{\Pi} e^{2 \varsigma\xi_2}
\vert \nabla v\vert^2 \ d\xi}.
\end{array}
\end{equation*}
Thus, for $\varsigma<\frac{2}{3}$ (that we may suppose), the
bilinear form $a_\varsigma$ satisfies $$a(v, Mv) \ge \gamma
\Vert v \Vert^2 \quad\hbox{ for every } v\in E_\varsigma,$$ with
$\gamma>0$. Then, by virtue of Lemma~\ref{l2} there is a unique
solution $\widetilde{\widetilde{X}}_{(I)}$ in $E_\varsigma$ of
the variational equation
$$a_\varsigma(\widetilde{\widetilde{X}}_{(I)},v)=L_\varsigma(v)
\quad \forall v \in E_\varsigma^0,$$ from which follows that
$\widetilde{\widetilde{X}}_{(I)}$ is a weak solution of Problem
(\ref{7!1}). Let us set, for simplicity of notation,
$Y=\widetilde{\widetilde{X}}_{(I)}$. From $Y\in E_\varsigma$ we
deduce that $Y$ decays exponentially fast in the Dirichlet
integral, i.e., for any $\delta>0$, there is a constant
$C_\delta$ such that $$\int\limits_{\Pi^\delta} \vert \nabla Y
\vert^2 \ d\xi \le C_{\delta}\, e^{-2 \varsigma\delta},$$ where
$\Pi^\delta=(-\frac{1}{2}, \frac{1}{2})\times (\delta, \infty)$.
Since $e^{\varsigma\xi_2}\frac{d \overline{Y}}{d\xi_2}\in
L_2(0,\infty)$ it follows that $\overline{Y}$ admits a limit as
$\xi_2\to + \infty$, which we denote $C_{(I)}(F)$. We have
$Y-\overline{Y}\in E_\varsigma^0$ and we easily show that
$Y-C_{(I)}(F)\in E_\varsigma^0$. Consequently, for any
$\delta>0$, there is a constant $C_\delta$ such that
$$\int\limits_{\Pi^\delta} \vert Y-C_{(I)(F)} \vert^2 \ d\xi \le
C_{\delta}\, e^{-2 \varsigma\delta}.$$ Then, using the local
regularizing properties of the Laplace operator and the Sobolev
imbedding theorem (see, for instance, \cite{So},~\cite{So1}), we
deduce, that $\forall \delta>0$, $\forall \alpha \in
\mathbb{N}^2$, $$ |\partial^\alpha
(Y-C_{(I)(F)})(\xi_1,\xi_2)|\leq C_{\delta,\alpha}\,
e^{-\varsigma\xi_2} $$ for any $(\xi_1,\xi_2)\in\Pi$ with
$\xi_2\ge\delta$, where $C_{\delta,\alpha}$ is another constant
depending only $\delta$ and $\alpha$. The proposition is proved.
\vskip20pt

{\it Proof of Proposition} \ref{p2}. Let $s\in
C^\infty(\mathbb{R})$ be such that $ s(\xi_2)=0 $ if $ \xi_2<1 $
and $s(\xi_2)=1$ if $\xi_2>2$, and let
$h(\xi_2)=\left(\frac{1}{6}\xi_2^3+\frac{1}{2}C(F)\xi_2^2\right)
s(\xi_2)$. Consider the problem
\begin{equation}\label{7!2p}
\left\{\begin{array}{l} \Delta_\xi
Z=X-h^{\prime\prime}(\xi_2)\quad\hbox{in} \ \Pi,\\
Z=0\quad\hbox{on} \ \Gamma,\qquad \frac{\partial Z}{\partial
\xi_1}=0\quad\hbox{as} \ \xi_1=\pm\frac{1}{2}.
\end{array}\right.
\end{equation}
Since $X-h^{\prime\prime}(\xi_2)=X-(\xi_2+C(F))$ for $\xi_2>2$,
according to \eqref{7'}, $$X-h^{\prime\prime}=0 \quad\hbox{as}\
\xi_2\to+\infty,$$ up to exponentially small terms. We then can
show as for Proposition~\ref{p1} that Problem \eqref{7!2p}
admits a solution $Z$ which has the asymptotics $$Z=C_{(II)}(F)
\quad\hbox{as}\ \xi_2\to+\infty,$$ up to exponentially small
terms, where $C_{(II)}(F)$ is a constant depending on the
function $F$. More precisely, denoting $$
C_{(II)}(F)=\lim_{\xi_2 \to +\infty }
\int\limits_{-\frac{1}{2}}^{\frac{1}{2}} Z(\xi_1,\xi_2)\
d\xi_1\qquad (\xi_2>0) $$ there exists a positive constant
$\varsigma>0$ such that, $\forall \delta>0$, $\forall \alpha \in
\mathbb{N}^2$, $$
|\partial^\alpha\left(Z(\xi_1,\xi_2)-C_{(II)}(F)\right)|\leq
C_{\delta,\alpha} e^{-\varsigma\xi_2} $$ for any
$(\xi_1,\xi_2)\in\Pi$ with $\xi_2>\delta$, where
$C_{\delta,\alpha}$ is a constant depending only $\delta$ and
$\alpha$. Setting $\widetilde{\widetilde{X}}_{(II)}=Z+h$, we
obtain the statement. \vskip20pt
\begin{remark}
Note that the existence of periodic in $(n-1)$ variables solutions
and their behavior at infinity in $n$-dimensional semi-space are
studied, for instance, in {\rm \cite{LP}}. See also {\rm
\cite{N}}.
\end{remark}
\section*{Acknowledgments.}
This paper has been partially written during the stay of Gregory
A. Chechkin in the Blaise Pascal University (Clermont-Ferrand,
France) in June -- July 2004. He wants to express deep thanks
for the hospitality, for the support and for perfect conditions
to work. The final version was completed when Gregory A.Chechkin
was visiting the Bashkir State Pedagogical University (Ufa,
Russia) in November 2006.



\begin{thebibliography}{B$_3$}


\bibitem{APV}
Achdou Y., Pironneau O., Valentin F. Effective Boundary
Conditions for Laminar Flows over Rough Boundaries, {\it J.
Comp. Phys.} {\bf 147} (1998), 187--218.


\bibitem{ABLS}
Amirat Y., Bresch D., Lemoine J., Simon J. Effect of Rugosity on
a Flow Governed by Navier-Stokes Equations. {\it Quarterly of
Applied Mathematics} {\bf LIX}: 4 (2001), 769--785.

\bibitem{AmCheGa}
Amirat Y., Chechkin G.A., Gadyl'shin R.R. Asymptotics of Simple
Eigenvalues and Eigenfunctions for the Laplace operator in a
Domain with Oscillating Boundary {\it Computat. Math. Math.
Phys.} {\bf 46}: 3 (2006), 97--110. Translated from {\it Zh.
Vychisl. Mat. i Mat. Fiz.} {\bf 46}: 1 (2006), 102--115.

\bibitem{AS}
Amirat Y., Simon J. Riblets and Drag Minimization. In: {\it
Optimization Methods in PDE's} (Eds. S. Cox and I. Lasiecka;
Contemporary Mathematics AMS) 1997, pp. 9--17.




\bibitem{BV}
Babu\v ska I., Vyborny R. Continuous Dependence of Eigenvalues
on the Domains. {\it Czech. Math. J.} {\bf 15} (1965), 169--178.



\bibitem{BMS}
Belyaev, A.G., Mikheev A.G., Shamaev A.S. Plane Wave diffraction
by a Rapidly Oscillating Surface. {\it Comput. Math. Math.
Phys.} {\bf 32}: 8 (1992), 1121--1133. Translated from {\it Zh.
Vychisl. Mat. i Mat. Fiz.} {\bf 32}: 8 (1992), 1258--1272.

\bibitem{BPC}
Belyaev A.G., Piatnitski A.L., Chechkin G.A. Asymp\-to\-tic
Behavior of a Solution to a Boundary-value Problem in a
Perforated Domain with Oscillating Boundary. {\it Siberian
Mathematical Journal} {\bf 39}: 4 (1998), 621--644. Translated
from {\it Sibirskii Matematicheskii Zhurnal} {\bf 39}: 4 (1998),
730--754.




\bibitem{Bor}
Borisov D.I. Asymptotics and Estimates for Eigenelements of
Laplacian with Frequent Nonperiodic Interchange of Boundary
Conditions {\it Izvestia: Mathematics.} {\bf 67}: 6 (2003),
1101--1148. Translated from {\it Izvestia RAN. Ser. Mat.} {\bf
67}: 6 (2003), 23--70.

\bibitem{BLS}
Bouchitte G., Lidouh A., Suquet P. Homog\'en\'eisation de
fronti\`ere pour la mod\'elisation du contact entre un corps
d\'eformable non lin\'eaire et un corps rigide. {\it C. R. Acad.
Sc. Paris. S\'er. I} {\bf 313} (1991), 967--972.

\bibitem{B}
{Brizzi R., Chalot J.P.} {Homog\'en\'eisation de fronti\`ere.}
{\it Ricerche di Mat.} {\bf 46}: 2 (1997), 341--387.


\bibitem{CUMN}
Chechkin G.A. Splitting of a Multiple Eigenvalue in a Problem on
Concentrated Masses {\it Russian Mathematical Surveys} {\bf 59}:
4 (2004), 790--791. 791 Translated from {\it UMN} {\bf 59}: 4
(2004), 205--206.


\bibitem{C}
Chechkin G.A. Asymptotic Expansion of Eigenvalues and
Eigenfunctions of an Elliptic Operator in a Domain with Many
``Light'' Concentrated Masses Situated on the Boundary.
Two--Dimensional Case. {\it Izvestia: Mathematics} {\bf 69}: 4
(2005), 805--846. Translated from {\it Izvestia RAN. Ser. Mat.}
{\bf 69}: 4 (2005), 161--204.

\bibitem{CCHE1}
Chechkin G.A., Chechkina T.P. On Homogenization Problems in
Do\-mains of the ``Infusorium'' Type. {\it Journal of
Mathematical Sciences} {\bf 120}: 3 (2004), 1470--1482.
Translated from {\it Trudy Seminara Imeni I.G. Petrovskogo} {\bf
23} (2003), 379--400.

\bibitem{CCHE2}
Chechkin G.A., Chechkina T.P. Homogenization Theorem for
Problems in Domains of the ``Infusorian'' Type with
Uncoordinated Structure. {\it Journal of Mathematical Sciences}
{\bf 123}: 5 (2004), 4363--4380. Translated from {\it Itogi
Nauki i Tekhniki. Sovrem. Mat. Prilozh. Tematicheskie obzory.
Differ. Uravn. Chast. Proizvod. (Progress in Science and
Technology. Series on Contemporary Mathematics and its
Applications. Thematic Surveys. Partial Differential Equations).
VINITI} {\bf 2} (2003), 139--154.

\bibitem{CC1}
Chechkin G.A., Cioranescu D. Vibration of a Thin Plate with a
``Rough'' Surface. In: {\it Nonlinear Partial Differential
Equations and their Applications. Coll\`ege de France Seminar.
Volume XIV. Studies in Mathematics and its Applications.}
Elsevier, Amsterdam - London - New York - Tokyo, 2002, pp.
147--169.

\bibitem{CFP}
Chechkin G.A., Friedman A., Piatnitski A.L. The Boundary--Value
Problem in Domains with Very Rapidly Oscillating Boundary. {\it
Journal of Math. Anal. and Applic.} {\bf 231}: 1 (1999),
213--234.

\bibitem{FHL}
Friedman A., Hu B., Liu Y. A boundary Value Problem for the
Poisson Equation with Multi--Scale Oscillating Boundary. {\it J.
Differential Equations} {\bf 137}: 1 (1997), 54--93.

\bibitem{G1}
Gadyl'shin R.R. Asymptotics of the Minimum Eigenvalue for a
Circle with Fast Oscillating Boundary Conditions. {\it C. R.
Acad. Sci., Paris, S\'er. I.} {\bf 323}: 3 (1996), 319--323.

\bibitem{G2}
Gadyl'shin R.R. Boundary-value Problem for the Laplacian with
Ra\-pid\-ly Oscillating Boundary Conditions. {\it Dokl. Math.}
{\bf 58}: 2 (1998), 293--296. Translated from {\it Dokl. Ross.
Akad. Nauk} {\bf 362}: 4 (1998), 456--459.

\bibitem{G3}
Gadyl'shin R.R. On the Eigenvalue Asymptotics for Periodically
Clam\-ped Membranes. {\it St. Petersbg. Math. J.} {\bf 10}: 1
(1999), 1--14. Translated from {\it Algebra Anal.} {\bf 10}: 1
(1998), 3--19.



\bibitem{Ga1}
Gaudiello A. Asymptotic Behavior of Non--Homogeneous Neumann
Problems in Domains with Oscillating Boundary. {\it Ricerche di
Math.} {\bf 43} (1994), 239--292.

\bibitem{I0}
Il'in A.M. A Boundary--Value Problem for an Elliptic Equation of
Second Order in a Domain with a Narrow Slit. I. The
Two-Dimensional Case (Russian). {\it Mat. Sb. (N.S.)} {\bf
99(141)}: 4 (1976), 514--537.

\bibitem{I1}
Il'in A.M. A Boundary--Value Problem for an Elliptic Equation of
Second Order in a Domain with a Narrow Slit. II. Domain with a
Small Opening (Russian). {\it Mat. Sb. (N.S.)} {\bf 103(145)}: 2
(1977), 265--284.


\bibitem{I3}
Il'in A.M. {\it Matching of Asymptotic Expansions of Solutions
of Boun\-da\-ry-va\-lue Problems}, AMS, Providence, 1992.

\bibitem{JwMa}
J\"ager W., Mikeli\'c A. On the Roughness--Induced Effective
Boundary Conditions for an Incompressible Viscous Flow. {\it J.
Differential Equations} {\bf 170} (2001): 96--122.

\bibitem{KPV}
{Kohler W., Papanicolaou G., Varadhan S.} Boundary and Interface
Problems in Regions with Very Rough Boundaries. In {\it Multiple
Scattering and Waves in Random Media} (Eds. Chow P.L., Kohler
W.E., Papanicolaou G.C.). North-Holland, Amsterdam, 1981, pp.
165--197.

\bibitem{LP}
Landis E.M., Panasenko G.P. A Theorem on the Asymptotics of
Solutions of Elliptic Equations with Coefficients Periodic in
all Variables Except One. {\it Soviet Math. Dokl.} {\bf 18}: 4
(1977), 1140--1143. Translated from {\it Doklady AN SSSR} {\bf
18}: 4 (1977), 1140--1143.

\bibitem{Lions}
{Lions J.L.}, {\it Some Methods in the Mathematical Analysis of
Systems and their Control}, Kexue Chubanshe (Science Press),
Beijing, Gordon and Breach Science Publishers, New York, 1981.

\bibitem{LS}
Lobo-Hidalgo M., S\'anchez-Palencia E. Sur certaines
propri\'et\'es spectrales des perturbations du domaine dans les
probl\`emes aux limites. {\it Comm. Partial Differential
Equations} {\bf 4} (1979), 1085--1098.

\bibitem{MK}
Marchenko V.A., Khruslov E.Ya. {\it Boundary--Value Problems in
Domains with Fine--Grained Boundaries.} (Russian) Naukova Dumka,
Kiev, 1974.


\bibitem{NM2}
Mel'nyk T.A., Nazarov S.A. Asymptotic Behavior of the Solution
of the Neumann Spectral Problem in a Domain of ``Tooth Comb''
Type. {\it J. Math. Sci.} (New York) {\bf 85}: 6 (1997),
2326--2346. Translated from {\it Trudy Seminara Imeni I.G.
Petrovskogo} {\bf 19} (1996), 138--173, 347.

\bibitem{NO}
Nazarov S.A., Olyushin M.V. Perturbation of the Eigenvalues of
the Neumann Problem due to the Variation of the Domain Boundary.
{\it St. Petersburg Math. J.} {\bf 5}: 2 (1994), 371--387.
Translated from {\it Algebra i Analiz} {\bf 5}: 2 (1993),
169--188.

\bibitem{N}
Nazarov S.A. Binomial Asymptotic Behavior of Solutions of
Spectral Problems with Singular Perturbations. {\it Math.
USSR-Sb.} {\bf 69}: 2 (1991), 307--340. Translated from {\it
Mat. Sb.} {\bf 181}: 3 (1990), 291--320.


\bibitem{NK}
Nevard J., Keller J.B. Homogenization of Rough Boundaries and
Interfaces. {\it SIAM J. Appl. Math.} {\bf 57}: 6 (1997),
1660--1686.

\bibitem{O}
Oleinik O.A., Shamaev A.S., Yosifian G.A. {\it Mathematical
Problems in Elasticity and Homogenization}. North-Holland,
Amsterdam, 1992.

\bibitem{Pe}
P\'erez M.E. On the Whispering Gallery Modes on Interfaces of
Membranes Composed of Two Materials with Very Different
Densities. {\it Mathematical Models and Methods in Applied
sciences} {\bf 13}:1 (2003), 75--98.



\bibitem{SP}
S\'anchez-Palencia E. {\it Homogenization Techniques for
Composite Media}. Springer-Verlag, Berlin - New York, 1987.

\bibitem{So}
Sobolev S.L. {\it Some Applications of Functional Analysis in
Mathematical Physics.} Translated from the Third Russian
Edition. Translations of Mathematical Monographs, 90. AMS,
Providence,
1991. 

\bibitem{So1}
Sobolev S.L. {\it Selected Problems in the Theory of Functional
Spaces and Generalized Functions.} (Russian) Nauka, Moscow,
1989.

\bibitem{Tartar}
Tartar L. {\it Homogenization, Compensated Compactness and
$H$-measures.} CBMS Conference, UC Santa Cruz, June 27--July 1,
1993. Unpublished.
\end{thebibliography}
\end{document}